\documentclass[11pt]{article}

\usepackage{amsmath, amssymb, amsthm, enumerate}
\usepackage{graphicx, caption, subcaption}
\usepackage{xcolor}
\usepackage{cancel} %
\usepackage[pagebackref]{hyperref}

\usepackage{xspace}

\usepackage{tikz}
\usetikzlibrary{arrows.meta,calc,decorations.markings}

\usepackage[margin=1.5in, marginparwidth=1.3in]{geometry}

\usepackage{titlesec}
\theoremstyle{plain}
\newtheorem{thm}{Theorem}[section]
\newtheorem{prop}[thm]{Proposition}
\newtheorem{lemma}[thm]{Lemma}

\newtheorem{claim}[thm]{Claim}

\theoremstyle{definition}

\newtheorem{rmkx}[thm]{Remark}
\newenvironment{rmk}
{\pushQED{\qed}\renewcommand{\qedsymbol}{$\triangle$}\rmkx}
{\popQED\endrmkx}

\renewcommand\qedsymbol{$\blacksquare$}

\newcommand{\C}{\mathbb{C}}

\newcommand{\sys}{\operatorname{sys}}
\newcommand{\genus}{\operatorname{genus}}
\newcommand{\length}{\operatorname{length}}

\newcommand{\area}{\operatorname{area}}

\newcounter{bencomments}

\title{Translation surfaces with logarithmically growing systole}  
\author{Benjamin Dozier \thanks{Department of Mathematics, Cornell University, \href{mailto:benjamin.dozier@cornell.edu}{\nolinkurl{benjamin.dozier@cornell.edu}}. Research supported in part by NSF Grant DMS-2247244 and the Simons Foundation.
}}

\begin{document}
\maketitle

\begin{abstract}
 We exhibit a sequence of translation surfaces $T_1,T_2,\ldots$ with $T_g$ of genus $g$ and area $g$, such that the systole grows logarithmically with $g$.  They are produced by starting with a fixed translation surface, taking covers with the same topology as those discovered by Brooks, and then adding handles by gluing along slits.    
\end{abstract}

\section{Introduction}
\label{sec:intro}

A \emph{translation surface} is a compact Riemann surface \(X\), equipped with a holomorphic \(1\)-form \(\omega\) that is not identically zero.  This gives rise to a flat cone metric on the surface, with cone points at the zeros of \(\omega\), and trivial holonomy.  The \emph{systole} \(\sys (X,\omega)\) is the length of the shortest non-contractible loop on \(X\).  In this note, we study how large the systole can be, asymptotically in the high genus regime.  This parallels analogous work for general Riemannian metrics, hyperbolic metrics, and graphs.  Our main result reveals similar behavior to that seen for hyperbolic surfaces and regular graphs, and the short proof goes by transporting a result about hyperbolic surfaces to translation surfaces.  To make this similarity more evident, we choose to normalize a translation surface of genus \(g\) to have area \(g\).

\begin{thm}\label{thm:all-genera}
  There exists a constant \(C>0\) and translation surfaces \(T_1,T_2,\ldots\) with \(\genus(T_g) = g= \area(T_g)\), such that for all $g$,
  \begin{align*}
    \sys (T_g) \ge C \log g .  
  \end{align*}
\end{thm}

We could rephrase the above as giving surfaces with systolic ratio at least order \((\log g)^2/g\).
The theorem resolves the relevant lower bound part of \cite[Conjecture 1.6]{chmw24}.  

In the other direction, it was previously known that there exists a constant \(C'\) such that
\begin{align}
  \label{eq:lower}
  \sys S \le C' \cdot \log (\genus(S))
\end{align}
for any translation surface $S$ with genus \(g\ge 2\) and area \(g\).  Thus the above theorem is optimal up to constant factors.  
The bound \eqref{eq:lower} follows (approximating the translation surface metric by Riemannian metrics) from Gromov's general systolic bound \cite [Corollary 6.4 D'']{gromov83}\footnote{there appears to be a typo in that statement -- the log factor should be raised to a power equal to the dimension; see \cite[ Theorem 2.2]{ks05} for a corrected version in the form we use here}.

\begin{rmk}
  The surfaces that we will produce in the proof have number of cone points growing linearly with the genus.  It is not possible to choose surfaces satisfying the conclusion of the theorem that lie in minimal strata (i.e. the metric has a single cone point, which has large cone angle).  In fact there exists a constant $C$ such that $\sys(S)\le C$ for any $S$ in a minimal stratum with area equal to genus.  To see this, consider growing a ball centered at the singular point.  If it is embedded, its area is order $gr^2$, where $r$ is its radius.  For $r$ large enough (independent of $g$), this area becomes larger than the area of $S$, hence the ball cannot be embedded, and from this one can produce a non-contractible loop of length at most $2r$.  
\end{rmk}

\paragraph{Previous work.}

Translation surface systoles, in our sense, have been studied in several papers.  Judge and Parlier \cite{jp19} study systoles of genus \(2\) translation surfaces.  In \cite{chmw24}, a general study of the ``systolic landscape'' is undertaken; they prove some general bounds on systole, and further investigate the genus \(2\) case.  
Buser-Makover-Muetzel \cite{bmm26} recently gave a non-trivial asymptotic lower bound on the maximum systole; in our normalization it translates to order \(\sqrt{\frac{\log g}{\log \log g}}\).  Their construction uses graphs of high girth.  They also suggest that the best bound should be of order $\log g$.  

Boissy-Geninska \cite{bg21} study analogous questions for the \emph{relative} systole (i.e. the length of the shortest saddle connection).  

Analogous bounds for graphs were proved by Erdős-Sachs \cite{es63}: they exhibit regular graphs (of fixed degree) with girth (length of shortest cycle) growing logarithmically in the number of vertices.  Hyperbolic surfaces with systole growing logarithmically in genus were found by Brooks \cite{brooks88}.  These surfaces were also studied by Buser-Sarnak \cite{bs94}.

\paragraph{Summary of proof.}

Brooks exhibited a hyperbolic surface with a sequence of covers that have logarithmic growth of systole as a function of genus.  We fix a translation surface and use covers with the same topology.  There is a Bilipschitz map between the base translation surface and base hyperbolic surface, and this lifts to (uniformly) Bilipschitz maps between the analogous covering surfaces.  This produces a sequence of translation surfaces with logarithmic systole growth (Proposition \ref{prop:seq}).  To get surfaces in every genus, we find many well-separated slits on the cover surfaces we've produced, and glue them in pairs.  

\paragraph{Acknowledgements. }  I thank Mingkun Liu and Bram Petri for useful discussions.  

\paragraph{AI usage statement. }  The (rather simple) proof strategy for Proposition \ref{prop:seq} was suggested by GPT-6 Astra while doing a literature search for another project. GPT-6 was also used to produce the figure, based on a hand-drawn diagram. The details of the proofs (including Lemma \ref{lem:bilip}, Lemma \ref{lem:discs-hyp}, and Claim \ref{claim:discs}) were then developed by me without further AI use.   The writing of the paper was done by me.  Once complete, a draft was run through GPT-6, which found several minor mistakes that were then corrected manually.

\section{Sequence of surfaces of growing genus}
\label{sec:proof}

In this section we show a somewhat weaker result: the existence of a sequence of surfaces with genera tending to infinity and logarithmically growing systole (rather than surfaces in every genus).   We will then modify the resulting surfaces to get all genera in the next section.  

\begin{prop}\label{prop:seq}
  There exists a constant \(C>0\) and translation surfaces \(S_n\) with \(\genus(S_n) = \area(S_n)\to\infty \), such that for all $n$,
  \begin{align*}
    \sys (S_n) \ge C \cdot \log \left(\genus(S_n) \right).  
  \end{align*}
\end{prop}

Since we will deduce the desired result from a systole bound for certain hyperbolic surfaces, we will need to compare translation surface metrics with hyperbolic metrics.  
The metric on a translation surface is not a Riemannian metric (except in genus \(1\)), but it does give a way of determining lengths of piecewise smooth paths, and thus also a metric in the sense of metric spaces (the distance between two points is the infimum over the lengths of paths connecting the two points).  

\begin{lemma}\label{lem:bilip}
  Let $(X,\omega)$ be a translation surface of genus $g$, and let $Y$ be a compact hyperbolic surface of genus $g$.  Then there is some $L$ for which there exists an $L$-Bilipschitz homeomorphism $f: X \to Y$.   
\end{lemma}

\begin{proof}
  The challenge is dealing with the cone points of the metric on $X$, which occur at zeros of $\omega$.  Near a point $p\in X$ where $\omega$ has a zero of order $k$, we can find a holomorphic coordinate chart to the unit ball $\Delta \subset \C$ for which $\omega$ corresponds to the standard form $z^kdz$.  We use polar coordinates on $\Delta$, writing $z=re^{i\theta}$.  We begin by defining a map $g_p: \Delta \to \Delta$ by $(r,\theta) \mapsto (r^{k+1},\theta)$.  Note that $g_p$ is a homeomorphism; it is a diffeomorphism away from the origin, but not at the origin.  In polar coordinates, the standard Euclidean metric is $\rho_{Euc}:=dr^2 + r^2d\theta^2$.  We compute the pull-back of $\rho_{Euc}$:    
  \begin{align*}
    g_p^* \rho_{Euc} &= d(r^{k+1})^2 + (r^{k+1})^2 d\theta^2 \\
    &= (k+1)^2 r^{2k} dr^2 + r^{2(k+1)}d\theta^2.  
  \end{align*}

  To see how $g_p$ distorts lengths, we want to compare this to the metric on $\Delta$ determined by $z^kdz$, which in polar coordinates is $\rho_{k}:= r^{2k} (dr^2 + r^2 d\theta^2)$.  Thus, using the above, we get
  \begin{align*}
    \rho_{k} \le g_p^* \rho_{Euc} \le (k+1)^2  \rho_{k},
  \end{align*}
  hence $g_p$ distorts lengths by at most a factor of $k+1$.  Now we post-compose $g_p$ with any smooth embedding  $i: \Delta \to Y$.  Since $i$ is smooth, by compactness, on any strictly smaller disc in $\Delta$, it distorts lengths by a bounded amount.  Combining with the above, we get that $i \circ g_p: \Delta \to Y$ distorts lengths by a bounded amount on any strictly smaller disc in $\Delta$.

  Now we can do this process simultaneously for the finitely many zeros $p$ of $\omega$, and composing with the coordinate maps on $X$, we get for each \(p \), a map \(f_p\) from an open neighborhood \(B_p\) of \(p\) into \(Y\).  To stitch these maps together, we make the \(B_p\) smaller so that they are disjoint in $X$, and shrink further to arrange that the images \(f_p(B_p)\) are disjoint in $Y$.  We then define $f$ to agree with $f_p$ on each $B_p$; on the complement of these discs, there is a smooth extension which restricts to a diffeomorphism $X-\bigcup_p B_p \to Y- \bigcup_p f_p\left(B_p \right)$.  By the preceding discussion, the map $f:X\to Y$ distorts lengths by a bounded amount on $\bigcup_pB_p$, and for the complement, both metrics on the source and target are Riemannian so the diffeomorphism property gives bounded distortion. 
  Thus $f$ is $L$-Bilipschitz for some $L$.  
 
\end{proof}

We can also think of the map in the proof above in flat geometric terms.  A zero \(p\) of \(\omega\) of order \(k\) has a neighborhood that is isometric to a neighborhood of the origin in the space obtained by gluing \(k+1\) copies of the Euclidean plane along rays through the origin.  We can then define a map on this space that ``reduces angles by a factor of \(1/(k+1)\)''.  It is intuitively reasonable that this map is (\(k+1\))-Bilipschitz.  

\begin{rmk}
  Note that the conclusion of the lemma need not hold if we replace \((X,\omega)\) by an arbitrary metric, in the sense of metric spaces.   In fact, consider the distance function \(d(x,y)\) associated to some Riemannian metric, and then replace it by \(d'(x,y):=\sqrt{d(x,y)}\).  This \(d'\) is also a metric (and gives the usual topology on the surface), but its Hausdorff dimension is \(4\); since \(d\) has Hausdorff dimension \(2\), the two are not Bilipschitz.  
\end{rmk}

\begin{proof}[Proof of Proposition \ref{prop:seq}]
  By \cite[Lemma 3.1]{brooks88} %
  there exists a hyperbolic surface \(Y\), covers \(p_n:Y_{n} \to Y\) where \(Y_{n}\) has genus \(g_n\to\infty\), and a constant \(C'>0\) such that  for all \(n\),
  \[\sys(Y_{n})\ge C' \log \left(\genus(Y_n)\right) .\]

  Now pick any translation surface \(S=(X,\omega)\) with \(\genus(S) = \genus(Y)\).  Let \(f: S \to Y\) be the \(L\)-Bilipschitz homeomorphism given by Lemma \ref{lem:bilip}.  Consider the covers \(q_n: S_n \to S\) with the same topology as \(p_n\);
  pull-back of \(\omega\) yields a translation surface structure on \(S_n\).   The map \(f\) lifts to an \(L\)-Bilipschitz homeomorphism \(S_n \to Y_n\).  It then follows that
\begin{align*}
  \sys(S_n) &\ge (1/L) \sys (Y_n) \ge (C'/L) \log \left(\genus(Y_n)\right) = (C'/L)  \log \left(\genus(S_n) \right).
\end{align*}
Now, since we are looking at covers, \(\area(S_n)\) is a linear function of \(\genus(S_n)\), so if we rescale \(S_n\) to have area equal to \(\genus(S_n)\), we get the same bound as in the previous inequality, with a different constant.  

\end{proof}

\section{All genera}
\label{sec:all-genera}

To prove Theorem \ref{thm:all-genera}, our task is to modify the surfaces produced in the proof of Proposition \ref{prop:seq}, increasing their genus without decreasing the systole too much.  We follow the ``gluing on handles'' strategy of \cite{ks25}, with some modifications needed for the translation surface context.  We will find many well-separated slits on our surfaces and then glue them in pairs to increase the genus.   First we find many well-separated discs on any hyperbolic surface with large systole; these will later be transported to the translation surfaces via the uniformly Bilipschitz maps, and the required slits will be found within them.

\begin{lemma}\label{lem:discs-hyp}
  Given $s,\epsilon>0$, there exists a constant $K>0$ such that on any hyperbolic surface $X$ of genus $g$ satisfying $\sys(X)\ge s \log g$, we can find $g^{1-\epsilon}$ disjoint embedded balls of radius $K\log g$.  
\end{lemma}

\begin{proof}
  We follow \cite{ks25}; we do not need as precise control as they do.   Let $B(r)$ be the ball of radius $r$ in the hyperbolic plane.  We have
\begin{align*}
  \area(B(r)) = 4\pi \sinh^2(r/2), \\
  \area(X) = 4\pi (g-1).
\end{align*}
Suppose that $r>0$ is such that
\begin{align}
  \label{eq:half-sys}
  r<\frac{s \log g}{2}.
\end{align}
  By the systole condition, any ball of radius $r$ on $X$ is embedded.  
Now consider a maximal collection of disjoint radius $r$ balls $B_1,\ldots,B_n$ on $X$, i.e. no other ball of this radius is disjoint from these.  Note that the balls $2B_i$ of double the radius and same centers must cover all of $X$, by the maximality condition.  It follows that the sum of the areas of the $2B_i$ is at least the area of $X$, hence using the area formulas above, we get 
\begin{align*}
  n\cdot  4\pi \sinh^2(r) \ge 4\pi(g-1).  
\end{align*}
Using the approximation $\sinh x \le e^x/2$, from the above we get
\begin{align}
  \label{eq:area}
  n\ge \frac{g-1}{\sinh^2r} \ge \frac{g-1}{e^{2r}/4} \ge \frac{2 g}{e^{2r}}.
\end{align}
Now for $K>0$ small in terms of $\epsilon,s$, setting $r=K\log g$, we can ensure both that \eqref{eq:half-sys} holds, and that $2r\le \epsilon \log g + \log 2$, which implies $2g^\epsilon \ge e^{2r}$, and hence $\frac{2g}{e^{2r}} \ge g^{1-\epsilon}.$  So from \eqref{eq:area}, we get $n\ge g^{1-\epsilon}$, and hence a maximal collection of radius $r$ balls has the required properties.  

\end{proof}

\begin{proof}[Proof of Theorem \ref{thm:all-genera}]

  We consider the surfaces in \cite[Lemma 3.1]{brooks88} constructed via congruence covers of compact arithmetic hyperbolic surfaces. We need some specific information about their genera.  Brooks finds a surface $Y$ and a sequence of covers $\pi_p : Y_p \to Y$ for each sufficiently large prime $p$, where $Y_p$ has genus  $g_p:=(p^3-p)\nu +1$, for some fixed positive integer $\nu$, and
  \begin{align}
    \label{eq:sys-log}
    \sys Y_p \ge \log g_p.  %
  \end{align}
    As in the proof of Proposition \ref{prop:seq}, we can find an $L$ and for each $p$ a translation surface $S_p$ equipped with a  $L$-Bilipschitz homeomorphism $f_p:S_p \to Y_p$.

    \begin{figure}[]
  \centering
  \resizebox{0.7\textwidth}{!}{%
\begin{tikzpicture}[
  x=1cm, y=1cm,
  line cap=round, line join=round,
  font=\normalsize,
  surface/.style={draw=black, line width=.85pt},
  disk/.style={draw=diskblue, line width=.8pt},
  slit/.style={draw=slitred, line width=1pt},
  disk label/.style={text=diskblue, inner sep=2pt},
  slit label/.style={text=slitred, inner sep=2pt}
]
\definecolor{diskblue}{HTML}{2466AD}
\definecolor{slitred}{HTML}{C43E32}

\draw[surface]
  (11.65,5.85)
  .. controls (10.25,5.98) and (7.80,5.85) .. (6.30,5.85)
  .. controls (4.65,5.85) and (3.00,5.90) .. (1.95,5.40)
  .. controls (.65,4.80) and (.25,3.85) .. (.25,2.95)
  .. controls (.25,1.78) and (.82,.85) .. (2.05,.40)
  .. controls (3.18,-.01) and (4.00,.15) .. (5.18,.38)
  .. controls (6.37,.61) and (7.06,.62) .. (8.17,.36)
  .. controls (9.45,.05) and (10.63,.11) .. (11.65,.48);

\draw[surface]
  (1.08,3.03) .. controls (1.39,2.77) and (2.20,2.77) .. (2.57,3.03);
\draw[surface]
  (1.36,3.00) .. controls (1.57,3.39) and (2.13,3.39) .. (2.32,3.00);

\draw[disk]
  (3.90,5.51)
  .. controls (4.39,5.63) and (5.03,5.44) .. (5.27,5.08)
  .. controls (5.51,4.73) and (5.21,4.28) .. (5.04,3.97)
  .. controls (4.84,3.60) and (4.20,3.49) .. (3.77,3.59)
  .. controls (3.49,3.66) and (3.44,3.81) .. (3.17,3.83)
  .. controls (2.82,3.85) and (2.70,4.21) .. (2.84,4.53)
  .. controls (2.96,4.80) and (2.64,4.88) .. (2.92,5.14)
  .. controls (3.20,5.40) and (3.52,5.51) .. (3.90,5.51)
  -- cycle;
\node[disk label] at (5.60,5.24) {$D_i$};

\draw[disk]
  (3.89,2.70)
  .. controls (4.44,2.81) and (5.02,2.57) .. (5.20,2.21)
  .. controls (5.37,1.91) and (5.13,1.62) .. (5.30,1.20)
  .. controls (5.43,.85) and (4.94,.69) .. (4.43,.64)
  .. controls (3.91,.55) and (3.35,.69) .. (3.19,1.06)
  .. controls (3.09,1.36) and (2.65,1.41) .. (2.75,1.78)
  .. controls (2.86,2.15) and (3.28,2.59) .. (3.56,2.67)
  .. controls (3.67,2.71) and (3.78,2.71) .. (3.89,2.70)
  -- cycle;
\node[disk label] at (5.61,2.29) {$D_j$};

\draw[disk]
  (7.75,5.48)
  .. controls (8.20,5.57) and (8.52,5.29) .. (8.98,5.14)
  .. controls (9.49,4.96) and (9.55,4.60) .. (9.45,4.17)
  .. controls (9.37,3.78) and (8.85,3.59) .. (8.48,3.65)
  .. controls (8.11,3.71) and (7.69,3.88) .. (7.58,4.12)
  .. controls (7.45,4.38) and (7.68,4.50) .. (7.39,4.63)
  .. controls (6.97,4.81) and (7.32,5.49) .. (7.75,5.48)
  -- cycle;
\draw[disk]
  (8.69,2.88)
  .. controls (9.03,2.98) and (9.40,2.69) .. (9.70,2.48)
  .. controls (10.02,2.25) and (10.04,1.86) .. (10.09,1.48)
  .. controls (10.17,1.00) and (9.65,.80) .. (9.23,.80)
  .. controls (8.83,.80) and (8.37,.75) .. (8.29,1.05)
  .. controls (8.18,1.45) and (7.73,1.55) .. (7.78,1.89)
  .. controls (7.83,2.29) and (8.25,2.74) .. (8.69,2.88)
  -- cycle;

\draw[slit] (4.01,4.10) -- (4.01,4.96);
\node[slit label,above] at (4.01,5.00) {$s_i$};
\node[slit label,left=5pt]  at (4.01,4.53) {$A$};
\node[slit label,right=5pt] at (4.01,4.53) {$B$};

\draw[slit] (4.01,1.23) -- (4.01,2.09);
\node[slit label,above] at (4.01,2.13) {$s_j$};
\node[slit label,left=5pt]  at (4.01,1.66) {$B$};
\node[slit label,right=5pt] at (4.01,1.66) {$A$};

\draw[slit] (8.39,4.11) -- (8.39,4.97);
\node[slit label,left=5pt]  at (8.39,4.54) {$C$};
\node[slit label,right=5pt] at (8.39,4.54) {$D$};

\draw[slit] (8.96,1.35) -- (8.96,2.21);
\node[slit label,left=5pt]  at (8.96,1.78) {$D$};
\node[slit label,right=5pt] at (8.96,1.78) {$C$};

\node at (11.00,3.10) {$\cdots$};
\end{tikzpicture}
  }
  \caption{Discs and slits produced in Claim \ref{claim:discs} and used in proof of Theorem \ref{thm:all-genera}.  The surface will be slit along each $s_i$, and then the sides of slits glued according to the labels.  Slitting along a pair and gluing increases the genus by one, and adds two new cone points of angle $4\pi$. }
  \label{fig:translation-systole-gluing}
\end{figure}
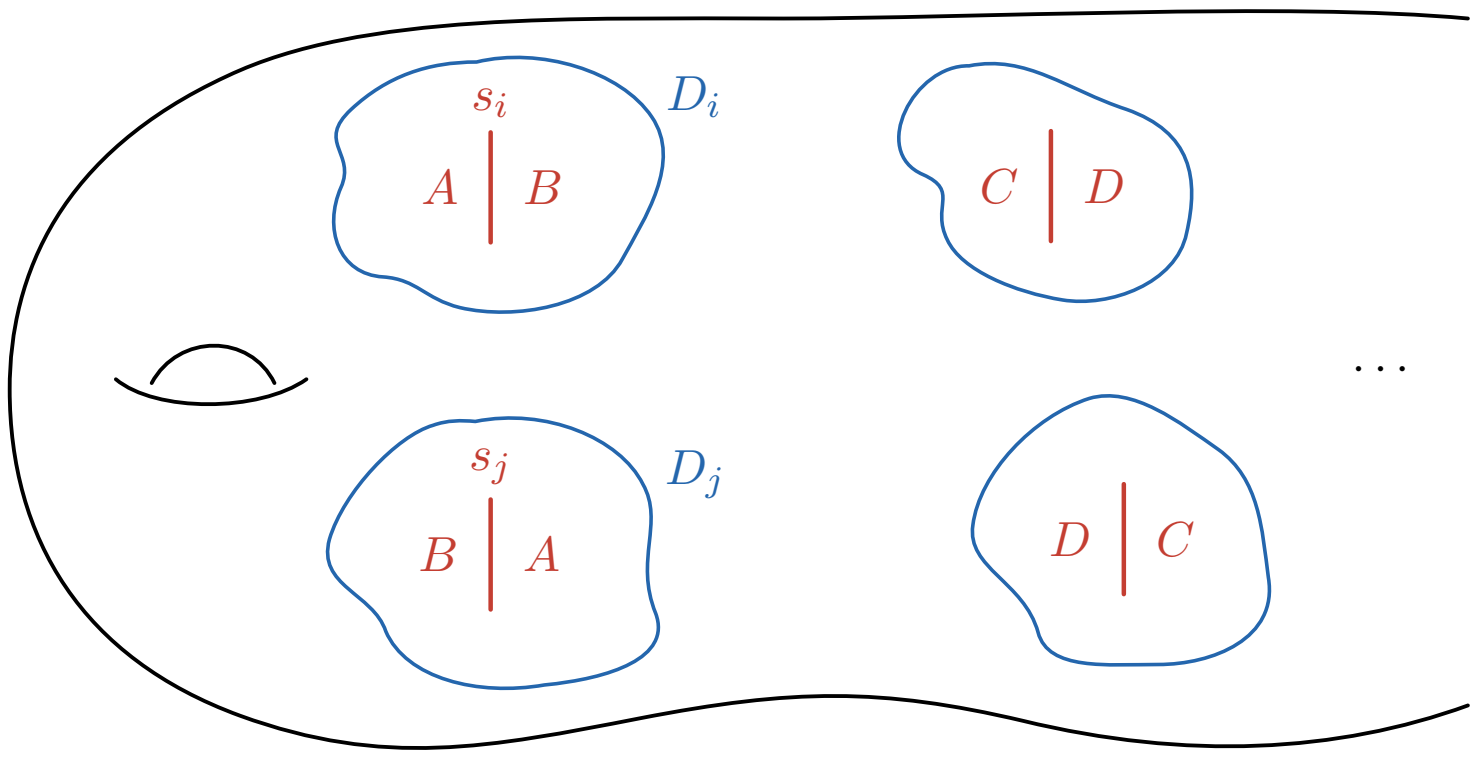

  \begin{claim} \label{claim:discs}
    For any $\epsilon>0$, there exists $c>0$ such that for each sufficiently large prime $p$, we can find $g_p^{1-\epsilon}$ disjoint topological discs $D_i$ on $S_p$ satisfying both:
    \begin{enumerate}[(i)]
    \item Each $D_i$ contains a geodesic segment $s_i$ of length $c\log g_p$ that does not contain a cone point.   The segments can be chosen to all be parallel to one another.  
    \item  $d(s_i,\partial D_i) \ge c \log g_p$ (here we are measuring the distance between the two sets with respect to the translation surface metric).    \label{item:bdy}
    \end{enumerate}
  \end{claim}

  \renewcommand\qedsymbol{$\square$}
  \begin{proof}
    Apply Lemma \ref{lem:discs-hyp} to the surfaces $Y_p$ (which by \eqref{eq:sys-log} satisfy a systole lower bound of the required form), giving a constant $K>0$ and $g_p^{1-\epsilon}$ disjoint embedded balls $D'_i$ on $Y_p$ of radius $r=K\log g_p$.  We take $D_i:=f_p^{-1}(D_i')$, which is a topological disk, since $D_i'$ is and $f_p$ is a homeomorphism.  Now consider the ball $\frac{1}{2}D_i'$ with the same center as $D_i'$ and half the radius (i.e. $r/2$).  Note that $f_p^{-1}\left(\frac{1}{2}D_i'\right)$ has distance at least $\frac{r}{2L}$ to $\partial D_i$.

    On the other hand, $f_p^{-1}\left(\frac{1}{2}D_i'\right)$ contains a ball $B_i$ of radius $\frac{r}{2L}$ in $S_p$.  In each $B_i$ we can find a geodesic segment $s_i$ of length $\frac{r}{3L}$ %
    avoiding cone points, and such that all the different $s_i$ are parallel (take a point near the center of each $B_i$ and not at a cone point, and shoot off a ray in some direction, avoiding the countably many bad directions that land at a cone point on some $B_i$).  Since $s_i$ is contained in $B_i$, it also has distance at least  $\frac{r}{2L}$ to $\partial D_i$.

    Taking $c$ to be $\frac{K}{3L}$ gives the desired result.  
  \end{proof}
\renewcommand\qedsymbol{$\blacksquare$}

  We will need to have some control over the gaps between the genera of this sequence of surfaces.  Let $p_1=2, p_2=3, p_3=5,\ldots$ be the sequence of primes.  We use Hoheisel’s bound \cite{hoheisel1930} on prime gaps, which implies  
  \begin{align}
    \label{eq:4}
    p_{n+1}-p_n  = O(p_n^{0.99997})
  \end{align}
  (considerably better exponents are known, but the above is enough for us).  It follows that
  \begin{align}
    \label{eq:g-gap}
    g_{p_{n+1}} - g_{p_{n}}  = O(p_{n+1}^3 - p_n^3) = O(p_n^{2.99997}) = O(g_{p_n}^{1-\epsilon}),
  \end{align}
  for some $\epsilon>0$.  

  Now combining this with Claim \ref{claim:discs} (applied with exponent $\epsilon/2$), we get there exists $c>0$ such that for sufficiently large $n$, we can find $M_{p_n}$ topological discs $D_i$ on $S_{p_n}$, with $M_{p_n} \ge 2 ( g_{p_{n+1}} - g_{p_{n}})$, containing geodesic segments $s_i$ of length $c\log g_{p_n}$ that are at least $c \log g_{p_n}$ away from the boundary of $D_i$.  We can then construct a translation surface  $T_g'$ for each $g\in [g_{p_n},g_{p_{n+1}}]$ as follows: repeatedly take a new pair of the segments $s_i,s_j$, slit the surface along that pair, and then identify each side of $s_i$ with a side of $s_j$ (the standard translation surface slit construction; see Figure \ref{fig:translation-systole-gluing}).  Doing this to each successive pair increases the genus by $1$, so by the condition on $M_{p_n}$, we exhaust the whole claimed interval of genera.

  Now we claim that for $C:=\min (1/L, c)$ and $g\in [g_{p_n},g_{p_{n+1}}]$, we have
  \begin{align}
    \label{eq:sys-bound}
    \sys(T_g') \ge C \log g_{p_n}.
  \end{align}
  In fact, consider a non-contractible curve $\gamma$ on $T_g'$.  To get the desired systole lower bound, it suffices to assume that $\gamma$ is \emph{geodesic} (we can always tighten $\gamma$ to either a concatenation of saddle connections, with turning angle at each point at least $\pi$, or a circumference curve of a cylinder; there can be many such geodesic representatives, but they all have the same length, which is the minimal length in the homotopy class).    Then $\gamma$ satisfies at least one of the following properties, from which we deduce the desired length bound:
  
  \begin{enumerate}[(i)]
  \item $\gamma$ is disjoint from all the %
    $s_i$ (which become pairs of saddle connections).  In this case, $\gamma$ can be thought of as a curve on $S_{p_n}$.  It is still a closed geodesic there, so its length must be at least $\sys(S_{p_n}) \ge \frac{1}{L}\sys(Y_{p_n}) \ge  C \log g_{p_n}$, using \eqref{eq:sys-log}.  %

  \item $\gamma$ is contained in the union of some pair of discs $D_i,D_j$ glued together along the segments $s_i,s_j$.  Any non-contractible curve in this union is homotopic to (a multiple of) that formed from the two saddle connections obtained from gluing $s_i,s_j$.  This curve has length $2c\log g_{p_n}$, and since it is geodesic, we get $\length(\gamma) \ge 2c\log g_{p_n}\ge C\log g_{p_n}$.  
    
  \item $\gamma$ intersects $s_i$ and $\partial D_i$ for some $i$.  Then by \eqref{item:bdy} of Claim \ref{claim:discs}, $\length(\gamma) \ge c\log g_{p_n} \ge C \log g_{p_n}$.  
  \end{enumerate}

We've now established \eqref{eq:sys-bound} and combining it with the bound for the genera gaps \eqref{eq:g-gap} we get, for each $g\in [g_{p_n},g_{p_{n+1}}]$,
\begin{align}
  \label{eq:sys-Tg'}
  \sys(T_g') &\ge C \log g_{p_n} \ge C \log \left(g_{p_{n+1}}-O(g_{p_n}^{1-\epsilon})\right) \ge 
               C \log g_{p_{n+1}} -O(1) \\
  & \ge C \log g - O(1). 
\end{align}
This gives surfaces in all sufficiently large genera $g$.  We now rescale $T_g'$ to produce surfaces $T_g$ of area $g$.  Since $g=g_{p_n} + O(g_{p_n}^{1-\epsilon})$ and $\area(T_g')=\area(S_{p_n})$ is a linear function of $g_{p_n}$, the scale factors we use are uniformly bounded (away from $0,\infty$), independent of $n$ and $g$.  Thus $\sys(T_g)$ differs from $\sys(T_g')$ by a uniformly bounded factor.  Combining this with \eqref{eq:sys-Tg'} gives
\begin{align*}
  \sys(T_g)  \ge C \log g - O(1)
\end{align*}
for some new constant $C$.  Reducing the constant $C$ further, to handle both low genera and the $O(1)$ term in the above, we get the desired result.  
  
\end{proof}

\begin{rmk}
  The Brooks construction does not give surfaces of all genera, which is why the gluing handles strategy used above is needed. The existence of a hyperbolic surface \(Y\), a constant \(C''\), and surfaces \(Y_g\) covering \(Y\)  of all sufficiently high genus \(g\) with
  \begin{align*}
    \sys(Y_g) \ge C'' \log g
  \end{align*}
  would yield a different proof of Theorem \ref{thm:all-genera} -- one could use these $Y_g$ in the proof of Proposition \ref{prop:seq} in place of the Brooks surfaces.

 Whether such a family $\{Y_g\}$ exists is apparently unknown. 
\end{rmk}

 \bibliographystyle{amsalpha}
   \bibliography{sources}%

@article {hoheisel1930,
    AUTHOR = {Hoheisel, Guido},
     TITLE = {Primzahlprobleme in der {Analysis}},
   JOURNAL = {Sitzungsber. Preuss. Akad. Wiss. Phys.-Math. Kl.},
      YEAR = {1930},
     PAGES = {580--588},
}

@article {ks25,
    AUTHOR = {Katz, Mikhail G. and Sabourau, St\'{e}phane},
     TITLE = {Logarithmic systolic growth for hyperbolic surfaces in every
              genus},
   JOURNAL = {Proc. Amer. Math. Soc.},
  FJOURNAL = {Proceedings of the American Mathematical Society},
    VOLUME = {153},
      YEAR = {2025},
    NUMBER = {1},
     PAGES = {325--330},
      ISSN = {0002-9939,1088-6826},
   MRCLASS = {53C45},
  MRNUMBER = {4840279},
MRREVIEWER = {Chady\ El Mir},
       DOI = {10.1090/proc/16985},
       URL = {https://doi.org/10.1090/proc/16985},
}

@article {brooks88,
    AUTHOR = {Brooks, Robert},
     TITLE = {Injectivity radius and low eigenvalues of hyperbolic
              manifolds},
   JOURNAL = {J. Reine Angew. Math.},
  FJOURNAL = {Journal f\"{u}r die Reine und Angewandte Mathematik. [Crelle's
              Journal]},
    VOLUME = {390},
      YEAR = {1988},
     PAGES = {117--129},
      ISSN = {0075-4102,1435-5345},
   MRCLASS = {58G25 (35P15 53C20)},
  MRNUMBER = {953681},
MRREVIEWER = {Carolyn\ Gordon},
       DOI = {10.1515/crll.1988.390.117},
       URL = {https://doi.org/10.1515/crll.1988.390.117},
}

@article {gromov83,
    AUTHOR = {Gromov, Mikhael},
     TITLE = {Filling {R}iemannian manifolds},
   JOURNAL = {J. Differential Geom.},
  FJOURNAL = {Journal of Differential Geometry},
    VOLUME = {18},
      YEAR = {1983},
    NUMBER = {1},
     PAGES = {1--147},
      ISSN = {0022-040X,1945-743X},
   MRCLASS = {53C20 (53C21 57R99)},
  MRNUMBER = {697984},
MRREVIEWER = {Yu.\ Burago},
       URL = {http://projecteuclid.org/euclid.jdg/1214509283},
}

@article {jp19,
    AUTHOR = {Judge, Chris and Parlier, Hugo},
     TITLE = {The maximum number of systoles for genus two {R}iemann
              surfaces with abelian differentials},
   JOURNAL = {Comment. Math. Helv.},
  FJOURNAL = {Commentarii Mathematici Helvetici. A Journal of the Swiss
              Mathematical Society},
    VOLUME = {94},
      YEAR = {2019},
    NUMBER = {2},
     PAGES = {399--437},
      ISSN = {0010-2571,1420-8946},
   MRCLASS = {32G15 (30F10 53C22)},
  MRNUMBER = {3942788},
MRREVIEWER = {Yunhui\ Wu},
       DOI = {10.4171/CMH/463},
       URL = {https://doi.org/10.4171/CMH/463},
}

@article {chmw24,
    AUTHOR = {Columbus, Tobias and Herrlich, Frank and Muetzel, Bjoern and
              Weitze-Schmith\"usen, Gabriela},
     TITLE = {Systolic geometry of translation surfaces},
   JOURNAL = {Exp. Math.},
  FJOURNAL = {Experimental Mathematics},
    VOLUME = {33},
      YEAR = {2024},
    NUMBER = {3},
     PAGES = {357--378},
      ISSN = {1058-6458,1944-950X},
   MRCLASS = {30F10 (32G15 53C22)},
  MRNUMBER = {4790974},
       DOI = {10.1080/10586458.2022.2108946},
       URL = {https://doi.org/10.1080/10586458.2022.2108946},
}

@article {bmm26,
    AUTHOR = {Buser, Peter and Makover, Eran and Muetzel, Bjoern},
     TITLE = {Translation surfaces with large systoles},
   JOURNAL = {J. Topol. Anal.},
  FJOURNAL = {Journal of Topology and Analysis},
    VOLUME = {18},
      YEAR = {2026},
    NUMBER = {5},
     PAGES = {1489--1504},
      ISSN = {1793-5253,1793-7167},
   MRCLASS = {30F10 (32G15 53C22)},
  MRNUMBER = {5051786},
       DOI = {10.1142/S1793525325500220},
       URL = {https://doi.org/10.1142/S1793525325500220},
}

@article {bg21,
    AUTHOR = {Boissy, Corentin and Geninska, Slavyana},
     TITLE = {Systoles in translation surfaces},
   JOURNAL = {Bull. Soc. Math. France},
  FJOURNAL = {Bulletin de la Soci\'et\'e{} Math\'ematique de France},
    VOLUME = {149},
      YEAR = {2021},
    NUMBER = {2},
     PAGES = {417--438},
      ISSN = {0037-9484,2102-622X},
   MRCLASS = {32G15 (30F30)},
  MRNUMBER = {4307431},
MRREVIEWER = {Fr\'ed\'eric\ Palesi},
       DOI = {10.24033/bsmf.2833},
       URL = {https://doi.org/10.24033/bsmf.2833},
}

@article {es63,
    AUTHOR = {Erdős, Paul and Sachs, Horst},
     TITLE = {Regul\"are {G}raphen gegebener {T}aillenweite mit minimaler
              {K}notenzahl},
   JOURNAL = {Wiss. Z. Martin-Luther-Univ. Halle-Wittenberg Math.-Natur.
              Reihe},
  FJOURNAL = {Wissenschaftliche Zeitschrift der Martin-Luther-Universit\"at
              Halle-Wittenberg. Mathematisch-Naturwissenschaftliche Reihe},
    VOLUME = {12},
      YEAR = {1963},
     PAGES = {251--257},
      ISSN = {0138-1504},
   MRCLASS = {55.10 (05.40)},
  MRNUMBER = {165515},
MRREVIEWER = {W.\ T.\ Tutte},
}

@article {bs94,
    AUTHOR = {Buser, P. and Sarnak, P.},
     TITLE = {On the period matrix of a {R}iemann surface of large genus},
      NOTE = {With an appendix by J. H. Conway and N. J. A. Sloane},
   JOURNAL = {Invent. Math.},
  FJOURNAL = {Inventiones Mathematicae},
    VOLUME = {117},
      YEAR = {1994},
    NUMBER = {1},
     PAGES = {27--56},
      ISSN = {0020-9910,1432-1297},
   MRCLASS = {22E40 (14H15 14H42 32G20)},
  MRNUMBER = {1269424},
MRREVIEWER = {Jos\'e\ M.\ Mu\~noz Porras},
       DOI = {10.1007/BF01232233},
       URL = {https://doi.org/10.1007/BF01232233},
}

@article {ks05,
    AUTHOR = {Katz, Mikhail G. and Sabourau, St\'ephane},
     TITLE = {Entropy of systolically extremal surfaces and asymptotic
              bounds},
   JOURNAL = {Ergodic Theory Dynam. Systems},
  FJOURNAL = {Ergodic Theory and Dynamical Systems},
    VOLUME = {25},
      YEAR = {2005},
    NUMBER = {4},
     PAGES = {1209--1220},
      ISSN = {0143-3857,1469-4417},
   MRCLASS = {53C23 (37C35 53C20)},
  MRNUMBER = {2158402},
MRREVIEWER = {Victor\ Bangert},
       DOI = {10.1017/S0143385704001014},
       URL = {https://doi.org/10.1017/S0143385704001014},
}

\end{document}